\documentclass[11pt]{article}
\usepackage[top=1in, bottom=1in, left=1in, right=1in]{geometry}
\usepackage{tikz}

\usepackage{amsmath,amsfonts,amssymb,amscd,amsthm,mathrsfs}
\usepackage{extarrows,textcomp}
\usepackage{graphicx,subfigure,epstopdf}
\usepackage{diagbox,multirow}
\usepackage{enumerate}
\usepackage{algorithm,algorithmic}
\usepackage{caption}
\usepackage{float}
\usepackage{lineno}
\usepackage{xcolor}
\usepackage{makeidx,hyperref}

\newtheorem{theorem}{Theorem}[section]

\newtheorem{example}{Example}[section]

\numberwithin{equation}{section}

\DeclareMathOperator*{\argmin}{arg\,min}

\newcommand{\dd}{\,{\rm d}}

\makeatletter
\newcommand*{\extendadd}{
  \mathbin{
    \mathpalette\extend@add{}
  }
}
\newcommand*{\extend@add}[2]{
  \ooalign{
    $\m@th#1\leftrightarrow$%
    \vphantom{$\m@th#1\updownarrow$}
    \cr
    \hfil$\m@th#1\updownarrow$\hfil
  }
}
\makeatother

\begin{document}

\title{Adaptive Learning on the Grids for  Elliptic Hemivariational Inequalities}

\author{Jianguo Huang ({\tt jghuang@sjtu.edu.cn})
\vspace{0.1in}\\
School of Mathematical Sciences, and MOE-LSC,\\
 Shanghai Jiao Tong University, Shanghai 200240, China
   \vspace{0.1in}\\
  Chunmei Wang ({\tt chunmei.wang@ttu.edu})
\vspace{0.1in}\\
Department of Mathematics \& Statistics, \\
Texas Tech University, 1108 Memorial Circle, Lubbock, TX 79409, USA
\vspace{0.1in}\\
Haoqin Wang ({\tt wanghaoqin@sjtu.edu.cn})
\vspace{0.1in}\\
School of Mathematical Sciences, and MOE-LSC,\\
 Shanghai Jiao Tong University, Shanghai 200240, China
}

\date{}
\maketitle

\begin{abstract}
This paper introduces a deep learning method for solving an elliptic hemivariational inequality (HVI). In this method, an expectation minimization problem is first formulated based on the variational principle of underlying HVI, which is solved by stochastic optimization algorithms using three different training strategies for updating network parameters. The method is applied to solve two practical problems in contact mechanics, one of which is a frictional bilateral contact problem and the other of which is a frictionless normal compliance contact problem. Numerical results show that the deep learning method is efficient in solving HVIs and the adaptive mesh-free multigrid algorithm can provide the most accurate solution among the three learning methods.
\end{abstract}


\section{Introduction}\label{sec: intro}
This paper considers an efficient numerical method for solving an elliptic hemivariational inequality (HVI), which was first introduced by Panagiotopoulos in the context of engineering in the early 1980s \cite{Panagiotopoulos1983}. After Panagiotopoulos' pioneering work, a lot of mathematical effort has been devoted to the study of HVIs (e.g.,  \cite{Panagiotopoulos1993,NaniewiczPanagiotopoulos1995,MotreanuPanagiotopoulos1999,MigorskiOchalSofonea2013,Han2020}), and HVIs have become a powerful mathematical tool in contact mechanics. In practice, the solution to elliptic HVIs is only available numerically. There have been various numerical methods to approximate the solution of elliptic HVIs (e.g., see~\cite{HaslingerMiettinenPanagiotopoulos1999,HanMigorskiSofonea2014,HanSofonea2019,FengHanHuang2019,WangQi2020}). After discretization, we often need to solve a non-convex and non-smooth optimization problem, which is  challenging and technical to attack. One of the most popular methods for this non-convex optimization is the iterative convexification approach~\cite{TzaferopoulosMistakidisBisbosPanagiotopoulos1995,MistakidisPanagiotopoulos1997}, where a sequence of convex problems approaching the original non-convex problem is constructed and solved. Based on this iterative convexification approach, many HVIs describing contact problems were solved, such as the frictional bilateral contact problem~\cite{BarboteuBartoszKalita2013,HanSofoneaBarboteu2017}, the frictionless normal compliance contact problem~\cite{BarboteuBartoszKalitaRamadan2014,HanSofoneaBarboteu2017}, and the frictionless unilateral contact problem~\cite{HanSofoneaBarboteu2017}. Except for the iterative convexification approach, other alternative choices are the proximal bundle method (e.g., see~\cite{Makela2002}), the bundle Newton method (e.g., see~\cite{BaniotopoulosHaslingerMoravkova2005}), and the primal-dual active-set algorithm~\cite{Kovtunenko2011}. Recently, Feng, Han and Huang~\cite{FengHanHuang2019} used the double bundle method (cf. \cite{JokiBagirovKarmitsaMakelaTaheri2018}) to solve the discrete non-convex and non-smooth problem arising from discretization of some HVIs. As a whole, the existing numerical approaches mentioned above are rather involved to implement in practice.

With the advance of deep learning techniques originated in computer science, considerable attention in applied mathematics and computational mathematics has been drawn to apply deep learning to scientific computing, especially for numerical methods to solve the partial differential equations (PDEs). Neural network-based numerical methods for solving PDEs date back to the 1990s~\cite{LeeKang1990,DissanayakePhan-Thien1994,LagarisLikasFotiadis1998} and achieve significant improvement recently  (e.g., see~\cite{McFallMahan2009,ShekariBeidokhtiMalek2009,RuddFerrari2015,SirignanoSpiliopoulos2018,HanJentzenE2018,KhooLuYing2019}). In those methods, deep neural networks (DNNs) are applied to parametrize PDE solutions and appropriate parameters are identified by minimizing an optimization problem formulated from the given PDE. The key to the success of neural networks-based methods is the universal approximation property of DNNs (cf. \cite{HornikStinchcombeWhite1989,Cybenko1989,Hornik1991,Barron1993,Kurkova1992,Yarotsky2017,Yarotsky2018,EMaWu2018,
ShenYangZhang2019,MontanelliDu2019,MontanelliYang2019,LuShenYangZhang2020,EWojtowytsch2020,SiegelXu2020}). Though deep learning has made great achievements in solving PDEs, especially in high-dimensional cases~\cite{HanJentzenE2018}, it is also challenging to obtain a highly accurate solution. Recently, several deep learning frameworks have been developed along this line. For instance, E and Yu proposed the Deep Ritz method for high-dimensional variational problems~\cite{EYu2018}. Karniadakis and his team designed several neural networks by incorporating physical information~\cite{RaissiPerdikarisKarniadakis2019,JagtapKawaguchiKarniadakis2020,JagtapKarniadakis2020}. Gu, Yang and Zhou developed the SelectNet to adaptively choose training samples in the learning propocess~\cite{GuYangZhou2020}. Liang et al. discussed how to design data-driven activation functions in~\cite{LiangLyuWangYang2021}. Bao and his team introduced weak adversarial networks to find the weak solution of PDEs in~\cite{ZangBaoYeZhou2020}. Chen et al. proposed Friedrichs learning strategies to solve various PDEs in a unified way~\cite{ChenHuangWangYang2020}.
Huang, Wang and Yang combined deep learning with traditional iteration to devise the Int-Deep for solving low-dimensional PDEs with a finite element accuracy~\cite{HuangWangYang2020}. Dong and Li combined the ideas of domain decomposition and extreme learning to form a new deep learning method for PDEs~\cite{DongLi2020}.

Though deep learning-based methods have been proposed to handle variational inequalities in~\cite{HuangWangYang2020}, to the best of our knowledge, there is no study on HVIs in the literature. In this paper, we propose a deep learning-based method to solve an elliptic HVI based on its equivalent variational form~\cite{Han2020} and compare the numerical performance of three different training strategies for updating parameters. In our method, the solution space of the HVI is parameterized via DNNs and an approximation is found by minimizing an unconstrained expectation minimization problem, which can be solved by stochastic gradient descent methods or its variants (cf. \cite{BottouCurtisNocedal2018}). In particular, the unconstrained expectation minimization problem is reformulated based on the variational principle of the HVI. Therefore, the resulting deep learning optimization problem has a clear physical meaning. In the meantime, we compare three algorithms for training networks: the classical stochastic algorithm with basic neural network parametrization, a blockwise training algorithm for multi-block neural networks, and an adaptive mesh-free multigrid algorithm. As applications, we employ the deep learning method to a frictional bilateral contact problem and a frictionless contact  problem with normal compliance. As we shall see in the numerical experiments, the deep learning method is efficient in solving HVIs and the adaptive mesh-free multigrid algorithm can provide the most accurate solution among the three learning methods discussed. Finally, it deserves to emphasize that our method is easy to be implemented in programming and suitable for practical applications.

 The rest of this paper is organized as follows. In Section~\ref{sec: HVI apps}, an elliptic HVI and its applications in contact mechanics are introduced briefly. In Section~\ref{sec: DL}, the deep learning method as well as three training algorithms for solving the HVI are presented in details. In Section~\ref{sec:experiments}, two numerical examples are provided to demonstrate the efficiency of the deep learning methods. Finally, we summarize our work with a short conclusion in Section \ref{sec:conclusion}.

\section{The elliptic hemivariational inequality with applications}
\label{sec: HVI apps}
 In this section, we first introduce some notations. For a real Banach space $X$ equipped with a norm $\| \cdot \|_X$, denote by $X^*$ its dual space equipped with a norm $\| \cdot \|_{X^*}$. The notation $\langle \cdot,\cdot \rangle$ stands for the dual pairing between $X^*$ and $X$. These notations apply to a Hilbert space $H$ naturally. Denote by $\mathcal{L}(X,Y)$ the space of all continuous linear operators from one normed linear space $X$ to another normed linear space $Y$. We also use the standard notations for Sobolev spaces and norms or seminorms (cf. \cite{Adams1975}).  
For any locally Lipschitz continuous functional $j$ on a Banach space $X$, denote by $j^0(u;w)$ the generalized (Clarke) directional derivative of $j$ at $u$ in the direction $w$ (cf. \cite{Clarke1975,Clarke1983}); i.e.,
\[
   j^0(u;w) = \lim_{v \to u}\sup_{t \downarrow 0}\frac{j(v+tw)-j(v)}{t}, \quad u\in X, w\in X.
\]
\subsection{The elliptic hemivariational inequality and its equivalent minimization problem}
\label{subsec: HVI}


In this paper, we consider an elliptic HVI on a spatial domain $\Omega$ in a finite-dimensional Euclidean space. For simplicity, we denote the boundary or part of the boundary of the domain as $\Gamma$. Let $H$ be a Hilbert space. The elliptic HVI can be described as follows.

\vspace{0.25cm}
\noindent
\textbf{Problem 1} Find $u \in H$ such that
\begin{equation}\label{HVI1}
  \langle Au,v\rangle + \int_{\Gamma}j^0(\gamma_j u;\gamma_j v) \dd s \ge \langle f, v \rangle, \quad v \in H,
\end{equation}
where $j: \Gamma \times \mathbb{R}^m \to \mathbb{R}$ is a locally Lipschitz continuous functional for some positive integer $m$, $\gamma_j \in \mathcal{L}(H, L^2(\Gamma; \mathbb{R}^m))$.
\vspace{0.25cm}
In the study of Problem 1, we need the following assumptions~\cite{Han2020}:
\vspace{-0.25cm}
\begin{align*}
	&(A1)\; A: H \to H^* \ \text{is Lipschitz continuous and strongly monotone; i.e., for a constant}\ m_A >0,  \\
	&\qquad \qquad \qquad \qquad \qquad
	\langle Av_1 - Av_2, v_1-v_2 \rangle \ge m_A \| v_1 - v_2 \|_H^2, \quad  v_1, v_2 \in H.
	\\
	&(A2)\; j(\cdot,\boldsymbol{z})\ \text{is measurable on}\ \Gamma, \ \boldsymbol{z}\in\mathbb{R}^m; \text{there exists}\ \boldsymbol{z}\in  L^2(\Gamma; \mathbb{R}^m)\ \text{such that}\ j(\cdot, \boldsymbol{z}(\cdot))\in L^1(\Gamma);
	\\
	&\qquad 
	\text{there exist non-negative constants}\  c_0, c_1,\ \text{and}\ \alpha_j\ \text{such that}
    \\
    &\qquad \qquad \qquad \qquad \qquad
    \| \partial j(\boldsymbol{z}) \|_{\mathbb{R}^{m}} \le c_0 + c_1\| \boldsymbol{z}\|_{\mathbb{R}^m}, \quad  \boldsymbol{z} \in \mathbb{R}^m,
    \\
    &\qquad \qquad \qquad \qquad \qquad
    j^0(\boldsymbol{z}_1; \boldsymbol{z}_2-\boldsymbol{z}_1) + j^0(\boldsymbol{z}_2;\boldsymbol{z}_1-\boldsymbol{z}_2) \le \alpha_j \| \boldsymbol{z}_1 -\boldsymbol{z}_2 \|_{\mathbb{R}^m}^2, \quad  \boldsymbol{z}_1,\boldsymbol{z}_2 \in \mathbb{R}^m.
    \\
    &(A3)\; f \in H^*.
    \\
    &(A4)\; \text{Denote by}\ c_{\Gamma}\ \text{an upper bound of the norm of the operator}\ \gamma_j.\
    \text{There holds}
    \\
    &\qquad \qquad \qquad \qquad \qquad
    \| \gamma_j v\|_{L^2(\Gamma; \mathbb R^m)}\le c_{\Gamma} \| v \|_H, \quad  v \in H.
\end{align*}

We recall an important result on the solution existence and uniqueness for Problem 1 (cf. \cite{HanSofoneaBarboteu2017}).

\begin{theorem}\label{thm1}
Assume (A1)-(A4) and $\alpha_{j}c_{\Gamma}^2 < m_A$ hold true. For any given $f\in H^*$, Problem 1 has a unique solution.
\end{theorem}

Next, assume $A\in \mathcal{L}(H, H^*)$ is symmetric, i.e.,
\[
  \langle Av_1, v_2 \rangle = \langle Av_2, v_1 \rangle, \quad  v_1, v_2 \in H.
\]
We turn to consider the following minimization problem:

\vspace{0.25cm}
\noindent\textbf{Problem 2}
Find $u \in H$ such that
\begin{equation}\label{HVI2}
  u = \argmin_{v\in H} E(v),
\end{equation}
where
\[
  E(v) = \frac{1}{2}\langle Av, v \rangle + \int_{\Gamma}j(\gamma_j v)\dd s - \langle f, v \rangle, \quad v \in H.
\]

\vspace{0.25cm}
As given in \cite{Han2020}, we have the following equivalence result.

\begin{theorem}\label{thm2}
Assume (A1)-(A4), $\alpha_{j}c_{\Gamma}^2 < m_A$, and $A\in \mathcal{L}(H, H^*)$ is symmetric. Then Problem 2 is equivalent to Problem 1.
\end{theorem}

Note that the functional $E(v)$ in~\eqref{HVI2} is reformulated from~\eqref{HVI1} based on the variational principle. Hence, it has clear physical meaning in mechanics.

\subsection{Some applications in contact mechanics}
\label{subsec: apps}
Let $\Omega\subset \mathbb{R}^2$ be the reference configuration of a linear elastic body. We assume that $\Omega$ is an open, bounded, and connected domain with Lipschitz continuous boundary $\Gamma = \partial \Omega$. The boundary is made of three disjoint and measurable parts: $\Gamma_D$, $\Gamma_T$, and $\Gamma_C$ such that $\text{meas}(\Gamma_D)>0$ and $\text{meas}(\Gamma_C)>0$, 
Denote by $\cdot$ and $|\cdot|$ the canonical inner product and the induced norm, respectively. For a vector field $\boldsymbol{v}\in \mathbb{R}^2$, we use $v_{\nu} =\boldsymbol{v} \cdot \boldsymbol{\nu}$ for its normal component and $\boldsymbol{v}_{\tau}=\boldsymbol{v} - v_{\nu}\boldsymbol{\nu}$ for its tangential component, where $\boldsymbol{\nu}$ is the unit outward normal vector to $\Gamma$.
The linearized strain tensor associated with a displacement field $\boldsymbol{u}: \Omega \to \mathbb{R}^2$ is denoted by $\boldsymbol{\varepsilon}(\boldsymbol{u})$ and the stress field is denoted by $\boldsymbol{\sigma}: \Omega \to \mathbb{S}^2$. In addition, we assume a volume force of density $\boldsymbol{f}_0 \in L^2(\Omega, \mathbb{R}^2)$ acting in $\Omega$. Besides, the body is assumed to be fixed on $\Gamma_D$, is subject to an action of the surface traction of density $\boldsymbol{f}_2 \in L^2(\Gamma_T, \mathbb{R}^2)$ on $\Gamma_T$, and is in contact on $\Gamma_C$. We   usually drop the spatial variable $\boldsymbol{x}$ for simplicity when the dependence is clear without any confusion. 

\vspace{0.25cm}
To discuss the contact problem, we introduce a space
\begin{equation*}
  V = \{ \boldsymbol{v} \in H^1(\Omega; \mathbb{R}^2)| \boldsymbol{v} = \boldsymbol{0} \ \text{a.e. on}\ \Gamma_D \}
\end{equation*}
equipped with the inner product
\begin{equation*}
 (\boldsymbol{u}, \boldsymbol{v})_V =  \int_{\Omega} \boldsymbol{\varepsilon}(\boldsymbol{u}) \cdot \boldsymbol{\varepsilon}(\boldsymbol{v}) \dd x, \quad \boldsymbol{u}, \boldsymbol{v} \in  V,
\end{equation*}
and the associated norm $\| \boldsymbol{v} \|_V = \sqrt{(\boldsymbol{v}, \boldsymbol{v})_V}$. Thanks to $\text{meas}(\Gamma_D)>0$ and Korn's inequality~\cite{Brenner2004}, we know
$V$ is a Hilbert space with the norm $\| \cdot \|_V$. Besides, we need to introduce another Hilbert space $Q = L^2(\Omega; \mathbb{S}^2)$, which is equipped with the inner product
\[
  (\boldsymbol{\sigma}, \boldsymbol{\tau})_Q =
  \int_{\Omega} {\sigma}_{ij}(\boldsymbol x){\tau}_{ij}(\boldsymbol x) \dd x.
\]
Here and below, we use the Einstein summation convention, which means the summation is implied for an index exactly appeared two times in a quantity.

\vspace{0.25cm}
First of all, we introduce an HVI to describe the frictional bilateral contact problem (cf. \cite{FengHanHuang2019,HanSofoneaBarboteu2017}).
We let
\[
  V_1=\{ \boldsymbol{v} \in V |  \boldsymbol{v}_{\nu} = 0 \ \text{on} \ \Gamma_C \},
\]
and define
\begin{align*}
	&\langle A \boldsymbol{u},  \boldsymbol{v} \rangle = (\boldsymbol{\mathcal{F}}(\boldsymbol{\varepsilon}(\boldsymbol{u})), \boldsymbol{\varepsilon}(\boldsymbol{v}))_Q;
	\\
	&\gamma_{j_{\tau}}: V_1 \to L^2(\Gamma_C;\mathbb{R}^2) \ \text{such that} \ \gamma_{j_{\tau}}(\boldsymbol{v})=\boldsymbol{v}_{\tau};
	\\
	&\int_{\Gamma}j^0(\gamma_j  \boldsymbol{u} ;\gamma_j  \boldsymbol{v}) \dd s = \int_{\Gamma_C} j_{\tau}^0(\gamma_{j_{\tau}} \boldsymbol{u} ;\gamma_{j_{\tau}} \boldsymbol{v}) \dd s;
	\\
	&\langle  \boldsymbol{f},  \boldsymbol{v}\rangle = \int_{\Omega} \boldsymbol{f}_0 \cdot \boldsymbol{v} \dd x + \int_{\Gamma_T}\boldsymbol{f}_2 \cdot \boldsymbol{v} \dd s, \  \boldsymbol{v} \in V;
\end{align*}
where the linear elasticity operator $\boldsymbol{\mathcal{F}} = (F_{ijkl})_{1 \le i,j,k,l \le 2}: \Omega \times \mathbb{S}^2 \to \mathbb{S}^2$ is symmetric, bounded, and satisfies the following property~\cite{Han2020}
\begin{equation}\label{LEO}
	\boldsymbol{\mathcal{F}}(\boldsymbol{\sigma}) \cdot \boldsymbol{\sigma} \ge m_{\boldsymbol{\mathcal{F}}} | \boldsymbol{\sigma} |^2, \quad m_{\boldsymbol{\mathcal{F}}}>0, \ \boldsymbol{\sigma} \in \mathbb{S}^2,
\end{equation}
$j_{\tau}: \Gamma_C \times \mathbb{R}^2 \to \mathbb{R}$ is locally Lipschitz on $\mathbb{R}^2$ for a.e. $\boldsymbol{x} \in \Gamma_C$ and satisfies the assumption $(A2)$ with constants $c_0$, $c_1$ and $\alpha_{j_{\tau}}$(cf. \cite{FengHanHuang2019,HanSofoneaBarboteu2017}).
Then the frictional bilateral contact problem can be described as a HVI in the form
\vspace{0.25cm}
\begin{equation}\label{FBP2}
  (\boldsymbol{\mathcal{F}}(\boldsymbol{\varepsilon}(\boldsymbol{u})), \boldsymbol{\varepsilon}(\boldsymbol{v}))_Q + \int_{\Gamma_C} j_{\tau}^0(\boldsymbol{u}_{\tau};\boldsymbol{v}_{\tau}) \dd s \ge \langle \boldsymbol{f}, \boldsymbol{v} \rangle, \quad \boldsymbol{v} \in V_1.
\end{equation}

Let $\boldsymbol{u}\in V_1$ and $\lambda_1$ be the smallest positive eigenvalue of the eigenvalue problem
 \[
    \int_{\Omega} \boldsymbol{\varepsilon}(\boldsymbol{u}) \cdot \boldsymbol{\varepsilon}(\boldsymbol{v}) \dd x = \lambda_1 \int_{\Gamma_C} \boldsymbol{u}_{\tau} \cdot \boldsymbol{v}_{\tau} \dd s, \quad \boldsymbol{v} \in V_1.
 \]
Then the assumption $(A4)$ is satisfied if $c_{\Gamma} \geq \sqrt{1/\lambda_1}$ holds true.

The assumption $(A1)$ is satisfied with $m_A = m_{\boldsymbol{\mathcal{F}}}$. Following Theorem~\ref{thm1} and Theorem~\ref{thm2}, if $\alpha_{j_{\tau}} < \lambda_1 m_{\boldsymbol{\mathcal{F}}}$, \eqref{FBP2} has a unique solution and is equivalent to the following optimization problem:

\vspace{0.25cm}
\begin{equation}\label{FBP3}
  \boldsymbol{u} = \argmin_{\boldsymbol{v} \in V_1} E(\boldsymbol{v}),
\end{equation}
where
\[
  E(\boldsymbol{v}) = \frac{1}{2}(\boldsymbol{\mathcal{F}}(\boldsymbol{\varepsilon}(\boldsymbol{v})), \boldsymbol{\varepsilon}(\boldsymbol{v}))_Q + \int_{\Gamma_C} j_{\tau}(\boldsymbol{v}_{\tau}) \dd s - \langle \boldsymbol{f}, \boldsymbol{v} \rangle.
\]

\vspace{0.25cm}
Next, we introduce an HVI to describe the frictionless normal compliance contact problem (cf. \cite{FengHanHuang2019,HanSofoneaBarboteu2017}).
Define
\begin{align*}
	&\langle A{\boldsymbol{u}}, {\boldsymbol{v}} \rangle = (\boldsymbol{\mathcal{F}}(\boldsymbol{\varepsilon}(\boldsymbol{u})), \boldsymbol{\varepsilon}(\boldsymbol{v}))_Q;
  	\\
  	&\gamma_{j_{\nu}}: V \to L^2(\Gamma_C) \ \text{such that} \ \gamma_{j_{\nu}}(\boldsymbol{v})= \boldsymbol{v}_{\nu};
  	\\
  	&\int_{\Gamma}j^0(\gamma_j \boldsymbol{u} ;\gamma_j \boldsymbol{v}) \dd s = \int_{\Gamma_C} j_{\nu}^0(\gamma_{j_{\nu}} \boldsymbol{u} ;\gamma_{j_{\nu}} \boldsymbol{v}) \dd s;
  	\\
	&\langle {\boldsymbol{f}}, \boldsymbol{v} \rangle = \int_{\Omega} \boldsymbol{f}_0 \cdot \boldsymbol{v} \dd x + \int_{\Gamma_T}\boldsymbol{f}_2 \cdot \boldsymbol{v} \dd s, \  \boldsymbol{v} \in V;
\end{align*}
where the linear elasticity operator $\boldsymbol{\mathcal{F}} = (F_{ijkl})_{1 \le i,j,k,l \le 2}: \Omega \times \mathbb{S}^2 \to \mathbb{S}^2$ is symmetric, bounded, and satisfies \eqref{LEO}(cf. \cite{Han2020}).
$j_{\nu}: \Gamma_C \times \mathbb{R} \to \mathbb{R}$ is locally Lipschitz on $\mathbb{R}$ for a.e. $\boldsymbol{x} \in \Gamma_C$ and satisfies the assumption $(A2)$ with constants $c_0$, $c_1$ and $\alpha_{j_{\nu}}$ (cf. \cite{FengHanHuang2019,HanSofoneaBarboteu2017}).

Therefore, the frictionless normal compliance contact problem reads
\vspace{0.25cm}
\begin{equation}\label{FNC2}
  (\boldsymbol{\mathcal{F}}(\boldsymbol{\varepsilon}(\boldsymbol{u})), \boldsymbol{\varepsilon}(\boldsymbol{v}))_Q + \int_{\Gamma_C} j_{\nu}^0( \boldsymbol{u}_{\nu}; \boldsymbol{v}_{\nu}) \dd s \ge \langle \boldsymbol{f}, \boldsymbol{v} \rangle, \quad \boldsymbol{v} \in V.
\end{equation}

\vspace{0.25cm}
Let $\boldsymbol{u}\in V$ and $\lambda_2$ be the smallest positive eigenvalue of the eigenvalue problem
 \[
   \int_{\Omega} \boldsymbol{\varepsilon}(\boldsymbol{u}) \cdot \boldsymbol{\varepsilon}(\boldsymbol{v}) \dd x = \lambda_2 \int_{\Gamma_C}  \boldsymbol{u}_{\nu} \cdot  \boldsymbol{v}_{\nu} \dd s, \quad  \boldsymbol{v} \in V.
\]
Then the assumption $(A4)$ is satisfied if $c_{\Gamma} \geq \sqrt{1/\lambda_2}$ holds true.
The assumption $(A1)$ is satisfied with $m_A = m_{\boldsymbol{\mathcal{F}}}$. Following Theorem~\ref{thm1} and Theorem~\ref{thm2}, if $\alpha_{j_{\nu}} < \lambda_2 m_{\boldsymbol{\mathcal{F}}}$, \eqref{FNC2} has a unique solution and is equivalent to the following optimization problem:

\vspace{0.25cm}
\begin{equation}\label{FNC3}
  \boldsymbol{u} = \argmin_{\boldsymbol{v} \in V} E(\boldsymbol{v}),
\end{equation}
where
\[
  E(\boldsymbol{v}) =\frac{1}{2}(\boldsymbol{\mathcal{F}}(\boldsymbol{\varepsilon}(\boldsymbol{v})), \boldsymbol{\varepsilon}(\boldsymbol{v}))_Q + \int_{\Gamma_C} j_{\nu}( \boldsymbol{v}_{\nu}) \dd s - \langle \boldsymbol{f}, \boldsymbol{v} \rangle.
\]

\section{Deep learning-based methods for HVIs}\label{sec: DL}
\subsection{The deep learning method}
The main idea of deep learning-based HVIs solvers is to treat DNNs as an efficient parametrization of the solution space of an HVI. The HVI solution is identified via seeking a DNN ${\phi}(\boldsymbol{x},\boldsymbol{\theta})$ with input $\boldsymbol{x}$ and parameters $\boldsymbol{\theta}$ that minimizes the variational minimization problem related to the HVI. From the discussion in Section~\ref{sec: HVI apps}, we know that Problem 1 is equivalent to Problem 2. This motivates the following problem.

\vspace{0.25cm}
\noindent\textbf{Problem 3}
Find $\boldsymbol{\theta}^*$ such that
\begin{equation}\label{HVI3}
  \boldsymbol{\theta}^* = \argmin_{\boldsymbol{\theta}} E({\phi}(\boldsymbol{x};\boldsymbol{\theta})),
\end{equation}
where
\begin{equation}\label{Loss}
  E({\phi}(\boldsymbol{x};\boldsymbol{\theta})) = \frac{1}{2}\langle A{\phi}(\boldsymbol{x};\boldsymbol{\theta}), {\phi}(\boldsymbol{x};\boldsymbol{\theta}) \rangle  + \int_{\Gamma}j(\gamma_j {\phi}(\boldsymbol{x};\boldsymbol{\theta})) \dd s - \langle f, {\phi}(\boldsymbol{x};\boldsymbol{\theta}) \rangle.
\end{equation}
In contact mechanics, the first and the third terms on the right-hand side of \eqref{Loss} usually can be formulated as the integrals. Then the objective function in \eqref{Loss} can be viewed as a sum of expectations of several random variables, which can be solved by stochastic gradient descent methods or its variants (cf. \cite{BottouCurtisNocedal2018}).  We refer to Subsection 4.1 for details along this line.

In this paper, we use two neural networks to approximate the solution of the HVI introduced in Section~\ref{sec: HVI apps}. The first one is the residual neural network (ResNet) proposed in \cite{HeZhangRenSun2016}.

Mathematically, the ResNet can be formulated as follow~\cite{EMaWang2019}:
\[
	\boldsymbol{h}_0=\boldsymbol{V}\boldsymbol{x},\
	\boldsymbol{g}_\ell=\sigma(\boldsymbol{W}_\ell\boldsymbol{h}_{\ell-1}+\boldsymbol{b}_{\ell}),\ 
    \boldsymbol{h}_\ell=\boldsymbol{h}_{\ell-1}+\boldsymbol{U}_\ell\boldsymbol{g}_\ell,\ \ell=1,2,\dots,L,\
    {\phi}(\boldsymbol{x};\boldsymbol{\theta})=\boldsymbol{a}^T\boldsymbol{h}_L,
\]
where $\boldsymbol{V}\in \mathbb{R}^{N_0\times d}$, $\boldsymbol{W}_\ell\in \mathbb{R}^{N_{\ell}\times N_{0}}$, $\boldsymbol{U}_\ell\in \mathbb{R}^{N_{0}\times N_{\ell}}$, $\boldsymbol{b}_\ell\in \mathbb{R}^{N_\ell}$ for $\ell=1,\dots,L$, $\boldsymbol{a}\in \mathbb{R}^{{N_0}\times m}$. $\sigma(x)$ is a non-linear activation function. For the purpose of simplicity, we consider $N_0=N_\ell=N$ and $\boldsymbol{U}_\ell$ is set as the identity matrix. Here, $L$ is  the depth of the ResNet, and $N$ is the width of the network, and $\boldsymbol{\theta} =\{ \boldsymbol{V}, \boldsymbol{a}, \boldsymbol{W}_{\ell}, \boldsymbol{b}_{\ell}: 1 \le \ell \le L\}$ denotes the set of all parameters in $\boldsymbol{\phi}$, which uniquely determines the neural network.

\begin{figure}[H]
\centering
\includegraphics[height = 0.2\textwidth]{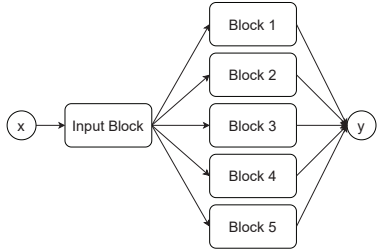}
\caption{The illustration of the network $B(\boldsymbol{x},\boldsymbol{\theta})$ with $P=5$.}
\label{fig: Block}
\end{figure}

The other one is a special neural network consisting of an input block and a few parallel blocks as visualized in Figure \ref{fig: Block}. 
Precisely speaking, the whole network is denoted by $B(\boldsymbol{x},\boldsymbol{\theta})$. The input block is a ResNet with width $N_{in}$ and depth $L_{in}$ denoted by $B_{in}(\boldsymbol{x},\boldsymbol{\theta}_{in})$. Following the input block are $P$ parallel blocks as independent ResNets with width $N_{p}$ and depth $L_{p}$ denoted by $B_{p}(\boldsymbol{x},\boldsymbol{\theta}_{p})$ for $1\le p \le P$. Then the network $B(\boldsymbol{x},\boldsymbol{\theta})$ can be formulated as
\begin{equation*}
\label{sturcture}
  B(\boldsymbol{x},\boldsymbol{\theta})= \sum_{p=1}^P B_{p}(B_{in}(\boldsymbol{x},\boldsymbol{\theta}_{in}),\boldsymbol{\theta}_{p})
\end{equation*}
In this network structure, different ResNets in parallel are trained with samples at different levels of discretization grids to obtain a PDE solution. The detailed training algorithm will be introduced later.

\subsection{Three strategies for updating network parameters}
Now let us introduce three different strategies to update the parameters of the networks.
First of all, we introduce a basic training algorithm to obtain a numerical solution to the HVI in the form of a ResNet ${\phi}(\boldsymbol{x};\boldsymbol{\theta})$ or the network ${B}(\boldsymbol{x};\boldsymbol{\theta})$. 
This algorithm is summarized in Algorithm~\ref{alg: basic} below.

\begin{algorithm}[H]
\caption{The basic training algorithm.}
\label{alg: basic}
\begin{algorithmic}
\REQUIRE The desired HVI.
\ENSURE The parameter set $\boldsymbol{\theta}$ in the solution ResNet ${\phi}(\boldsymbol{x};\boldsymbol{\theta})$ or the solution network ${B}(\boldsymbol{x};\boldsymbol{\theta})$.
\STATE Set $Epoch$ as the total iteration number, and assign $N$, $N_{\Gamma}$ as the sample sizes in the domain $\Omega$, the boundary $\Gamma$, respectively. 
\STATE Initialize ${\phi}(\boldsymbol{x};\boldsymbol{\theta})$ or ${B}(\boldsymbol{x};\boldsymbol{\theta})$ following the default random initialization of PyTorch.
\FOR{$k=1,\cdots,Epoch$}
\STATE Generate uniformly distributed samples $\{\boldsymbol{x}_i\}_{i=1}^{N}\subset \Omega$, $\{\boldsymbol{x}_i\}_{i=1}^{N_{\Gamma}}\subset\Gamma$.
\STATE Update $\boldsymbol{\theta}$ using the gradient of \eqref{Loss} evaluated at the generated samples.
\ENDFOR
\end{algorithmic}
\end{algorithm}

Next, we introduce a blockwise training algorithm without utilizing multigrid sampling to obtain the approximate solution to the HVI in the form of the network ${B}(\boldsymbol{x},\boldsymbol{\theta})$ introduced before. This algorithm can serve as a comparison to the adaptive multigrid training algorithm to be introduced later. It consists of two phases: the initialization phase and the refinement phase. In the initialization phase, we  train the network ${B}(\boldsymbol{x},\boldsymbol{\theta})$ by Algorithm~\ref{alg: basic} for $Epoch_{int}$ iterations. In the refinement phase, we train only one block ${B}_p(\boldsymbol{x},\boldsymbol{\theta}_p)$ together with the input block (i.e., the parameters in other blocks are fixed) in each iteration. We use $Epoch_{re}$ to denote the number of blockwise refinement steps and $Epoch_{b}$ to specify how many times a certain block is trained contiguously. Thus, the total training iterations in the refinement phase is $Epoch_{re}\times P \times Epoch_{b}$. The blockwise training algorithm is summarized in Algorithm~\ref{alg: blockwise} below.

\begin{algorithm}[H]
\caption{The blockwise training algorithm.}
\label{alg: blockwise}
\begin{algorithmic}
\REQUIRE The desired HVI.
\ENSURE The parameter set $\boldsymbol{\theta}$ in the solution network ${B}(\boldsymbol{x};\boldsymbol{\theta})$.
\STATE Set parameters $Epoch_{int}$, $Epoch_{re}$, and $Epoch_{b}$ as iteration numbers. Assign parameters $N$, $N_{\Gamma}$ for sample sizes in the domain $\Omega$, the boundary $\Gamma$, respectively. Let $P$ be the number of blocks. 
\STATE Initialize ${B}(\boldsymbol{x};\boldsymbol{\theta})$ following the default random initialization of PyTorch.
\FOR{$k=1,\cdots,Epoch_{int}$}
\STATE Generate uniformly distributed samples $\{\boldsymbol{x}_i\}_{i=1}^{N}\subset \Omega$, $\{\boldsymbol{x}_i\}_{i=1}^{N_\Gamma}\subset\Gamma$.
\STATE Update $\boldsymbol{\theta}$ using the gradient of \eqref{Loss} evaluated at the generated samples. 
\ENDFOR
\FOR{$j=1,\cdots,Epoch_{re}$}
\FOR{$p=1,\cdots,P$}
\FOR{$\ell=1,\cdots,Epoch_{b}$}

\STATE Generate uniformly distributed samples $\{\boldsymbol{x}_i\}_{i=1}^{N}\subset\Omega$, $\{\boldsymbol{x}_i\}_{i=1}^{N_{\Gamma}}\subset\Gamma$.
\STATE Update all the parameters of the $p$-th block and the input block using the gradient of \eqref{Loss} evaluated at the sampled grid points just above.
\ENDFOR
\ENDFOR
\ENDFOR
\end{algorithmic}
\end{algorithm}

 Since deep learning is challenging to obtain a very accurate solution to PDEs, we use some technical strategies for training the networks to improve the accuracy. It is well known that multigrid methods improve solution accuracy via computation on a hierarchy of discretization meshes. Multigrid computation traverses from fine grids to coarse grids and finally back to fine grids so as to efficiently improve the solution error at different levels. The special order of traverse in multigrid methods is motivated by the fact that traditional iterative methods tend to reduce high-frequency errors better than low-frequency errors. In contrast, deep learning optimization reduces low-frequency errors faster than high-frequency errors as discussed in~\cite{XuZhangLuoXiaoMa2020, CaoFangWuZhouGu2019}. Based on this observation, we propose a mesh-free multigrid method for training the parameters of the networks.
 A straightforward design of such a method would move the computation among different levels of grids in the opposite way of multigrid. However, the numerical performance is not satisfactory.
 To overcome this limitation, we adopt the adaptive strategy, that means, the computation  moves from the current grids to the grids associated with the smallest loss function evaluated at the corresponding grids. This idea leads to an adaptive multigrid method described below.

 Suppose there are $P$ levels of grids generated by uniform discretization with step sizes $H,2H$,$\cdots$, $2^{P-1}H$, where $H$ is the step size of the finest grids. Let $Epoch_{int}$, $Epoch_{re}$, and $Epoch_{b}$ be the same parameters as in the blockwise training algorithm. This algorithm utilizing multigrid sampling also consists of an initialization phase and a refinement phase. The initialization phase of the adaptive multigrid training is the same as the one of the blockwise training algorithm. In the refinement phase, we first evaluate a loss $E_p$ according to \eqref{Loss} using $N$ samples from $\Omega_p^H$, $N_{\Gamma}$ samples from $\Gamma_{p}^H$, where $\Omega^H_p$ is the grids in the domain with a step size $2^{P-1}H$, $\Gamma_{p}^H$ is  the grids on the boundary $\Gamma$ with a step size $2^{P-1}H$. If $E_k$ is the smallest one among $\{E_p\}_{p=1}^P$, we then train the $k$-th block and the input block with $N$ samples from $\Omega_p^H$, $N_{\Gamma}$ samples from $\Gamma_{p}^H$. Thus, the total training iterations in the refinement phase is $Epoch_{re}\times Epoch_{b}$. Such a training algorithm chooses the grids corresponding to the smallest loss as training samples to refine the network parameters and, hence, it is called an adaptive multigrid training algorithm. See Algorithm \ref{alg: adaptive}  below for a detailed description.

\begin{algorithm}[H]
\caption{The adaptive multigrid training.}
\label{alg: adaptive}
\begin{algorithmic}
\REQUIRE The desired HVI.
\ENSURE The parameter set $\theta$ in the solution network ${B}(\boldsymbol{x};\boldsymbol{\theta})$.
\STATE Set parameters $Epoch_{int}$, $Epoch_{b}$, and $Epoch_{re}$ as iteration numbers, parameters $N$ and $N_{\Gamma}$ as sample sizes in the domain $\Omega$ and the boundary $\Gamma$, respectively. Let $P$ be the number of blocks. 
\STATE Initialize ${B}(\boldsymbol{x};\boldsymbol{\theta})$ following the default random initialization of PyTorch.
\FOR{$k=1,\cdots,Epoch_{int}$}
\STATE Generate uniformly distributed samples $\{\boldsymbol{x}_i\}_{i=1}^{N}\subset\Omega$, $\{\boldsymbol{x}_i\}_{i=1}^{N_\Gamma}\subset\Gamma$.
\STATE Update $\boldsymbol{\theta}$ using the gradient of \eqref{Loss} evaluated at the generated samples.
\ENDFOR
\FOR{$j=1,\cdots,Epoch_{re}$}
\FOR{$p=1,\cdots,P$}
\STATE Generate uniformly distributed samples $\{\boldsymbol{x}_i\}_{i=1}^{N}\subset\Omega^H_p$, $\{\boldsymbol{x}_i\}_{i=1}^{N_{\Gamma}}\subset\Gamma_{p}^H$.
\\
Evaluate the loss $E_p$ according to \eqref{Loss} using the samples just above.
\ENDFOR
\STATE Let $p=\argmin_{1\leq p\leq P} E_p$.
\FOR{$\ell=1,\cdots,Epoch_{b}$}
\STATE Generate uniformly distributed samples $\{\boldsymbol{x}_i\}_{i=1}^{N}\subset\Omega^H_p$, $\{\boldsymbol{x}_i\}_{i=1}^{N_{\Gamma}}\subset\Gamma_{p}^H$.
\\
Update all the parameters of the $p$-th block and the input block using the gradient of \eqref{Loss} evaluated at the sampled grid points just above.
\ENDFOR
\ENDFOR
\end{algorithmic}
\end{algorithm}

\section{Numerical experiments}\label{sec:experiments}
In this section, we shall illustrate the performance of the deep learning method to a frictional bilateral contact problem and a frictionless normal compliance contact problem. Numerical comparisons of the deep learning method and the virtual element method (VEM) are provided. As we will see, the approximation accuracy in the form of DNNs is nearly the same as that of VEM in a fine meshsize, but the deep learning method is much easier to implement. In the meantime, we investigate the numerical performance in terms of different network structures and different training algorithms. As we will see, the adaptive mesh-free multigrid algorithm can provide a more accurate approximation to HVIs than other deep learning methods.

\subsection{The algorithm description of deep learning methods for contact problems}
From Section~\ref{sec: HVI apps}, we know the frictional bilateral contact problem~\eqref{FBP2} is equivalent to optimization problem \eqref{FBP3} and the frictionless normal compliance contact problem~\eqref{FNC2} is equivalent to \eqref{FNC3}. In order to numerically solve these problems by deep learning method, we need to parametrize the solution space of HVIs as introduced in Section~\ref{sec: DL} firstly. The solution spaces are $V_1$ for the frictional bilateral contact problem~\eqref{FBP3} and $V$ for the frictionless normal compliance contact problem~\eqref{FNC3}. To deal with the constraints in the admissible spaces $V_1$ or $V$, a natural way is the penalty method that penalizes the loss function with extra terms to enforce these constraints. However, tuning parameters in the penalty method may be tedious in practice. Therefore, we construct DNNs satisfying these constraints automatically as follows:
\[
  \boldsymbol{\phi}(\boldsymbol{x};\boldsymbol{\theta}) = \boldsymbol{b}(\boldsymbol{x})*\boldsymbol{\psi}(\boldsymbol{x};\boldsymbol{\theta}),
\]
where $\boldsymbol{b}(\boldsymbol{x})$ is a known smooth vector-valued function such that $\boldsymbol{\phi}=\boldsymbol{0}$ on the Dirichlet boundary $\Gamma_D$ and $\phi_{\nu}=0$ on the contact boundary $\Gamma_C$ for the frictional bilateral contact problem, or $\boldsymbol{\phi}=\boldsymbol{0}$ on the Dirichlet boundary $\Gamma_D$  for the frictionless normal compliance contact problem; $\boldsymbol{\psi}(\boldsymbol{x};\boldsymbol{\theta})$ is an arbitrary DNN in $\{ \boldsymbol{\psi}(\boldsymbol{x};\boldsymbol{\theta}) \}_{\boldsymbol{\theta}} \approx H^1(\Omega; \mathbb{R}^2)$; and ``$*$'' stands for the Hadamard product. Thus, \eqref{HVI3} can be formulated as follows

\begin{equation}\label{HVI4}
 \boldsymbol{\theta}^* = \argmin_{\boldsymbol{\theta}}E(\boldsymbol{b}(\boldsymbol{x})*\boldsymbol{\psi}(\boldsymbol{x}; \boldsymbol{\theta})),
\end{equation}
and $\boldsymbol{u}^{DL}(\boldsymbol{x};\boldsymbol{\theta}^*)=\boldsymbol{b}(\boldsymbol{x})*\boldsymbol{\psi}(\boldsymbol{x}; \boldsymbol{\theta}^*) $ is the approximate solution to the target HVI.

\vspace{0.25cm}
Next, we discuss the specific form $E(\boldsymbol{b}*\boldsymbol{\psi})$ in \eqref{HVI4} for the frictional bilateral contact problem~\eqref{FBP3} and the frictionless normal compliance contact problem~\eqref{FNC3}, respectively.

\vspace{0.25cm}
For the frictional bilateral contact problem~\eqref{FBP3},
\begin{align}\label{LossB}
& \quad E(\boldsymbol{b}(\boldsymbol{x})*\boldsymbol{\psi}(\boldsymbol{x}; \boldsymbol{\theta}))  \nonumber
\\
&= \frac{1}{2}(\boldsymbol{\mathcal{F}}(\boldsymbol{\varepsilon}(\boldsymbol{b}(\boldsymbol{x})*\boldsymbol{\psi}(\boldsymbol{x}; \boldsymbol{\theta})), \boldsymbol{\varepsilon}(\boldsymbol{b}(\boldsymbol{x})*\boldsymbol{\psi}(\boldsymbol{x}; \boldsymbol{\theta}))_Q
   + \int_{\Gamma_C} j_{\tau}((\boldsymbol{b}*\boldsymbol{\psi})_{\tau}(\boldsymbol{x}; \boldsymbol{\theta})) \dd s
   - \langle \boldsymbol{f}, \boldsymbol{b}(\boldsymbol{x})*\boldsymbol{\psi}(\boldsymbol{x}; \boldsymbol{\theta}) \rangle  \nonumber
   \\
   &=|\Omega|\mathbb{E}_{\boldsymbol{\xi}_1}\left[\frac{1}{2} \boldsymbol{\mathcal{F}}(\boldsymbol{\varepsilon}(\boldsymbol{b}(\boldsymbol{\xi_1})*\boldsymbol{\psi}(\boldsymbol{\xi}_1; \boldsymbol{\theta})))_{ij}  \boldsymbol{\varepsilon}(\boldsymbol{b}(\boldsymbol{\xi_1})*\boldsymbol{\psi}(\boldsymbol{\xi}_1; \boldsymbol{\theta}))_{ij} - \boldsymbol{f}_0(\boldsymbol{\xi}_1) \cdot \boldsymbol{b}(\boldsymbol{\xi_1})*\boldsymbol{\psi}(\boldsymbol{\xi}_1; \boldsymbol{\theta}) \right]   \nonumber
   \\ &\quad - |\Gamma_T|\mathbb{E}_{\boldsymbol{\xi}_2} \big[ \boldsymbol{f}_2(\boldsymbol{\xi}_2) \cdot \boldsymbol{b}(\boldsymbol{\xi_2})*\boldsymbol{\psi}(\boldsymbol{\xi}_2; \boldsymbol{\theta}) \big]
   +|\Gamma_C|\mathbb{E}_{\boldsymbol{\xi}_3}\big[ j_{\tau}((\boldsymbol{b}*\boldsymbol{\psi})_{\tau}(\boldsymbol{\xi}_3; \boldsymbol{\theta})) \big],
\end{align}
with $\boldsymbol{\xi}_1$, $\boldsymbol{\xi}_2$, and $\boldsymbol{\xi}_3$ being random vectors following the uniform distribution over $\Omega$, $\Gamma_T$, and $\Gamma_C$, respectively.
In practice, the minimization problem~\eqref{LossB} is solved by the stochastic gradient descent method \cite{BottouCurtisNocedal2018} or its variants (e.g. Adam\cite{KingmaBa2014}) by randomly sampling the integral domains in the loss function. In each iteration of the optimization algorithm, a stochastic loss function defined below is minimized instead of the original loss function in~\eqref{LossB}:
\begin{align}\label{FBP6}
  &\quad\hat{E}(\boldsymbol{b}(\boldsymbol{x})*\boldsymbol{\psi}(\boldsymbol{x}; \boldsymbol{\theta})) \nonumber
  \\
  &=\frac{|\Omega|}{N}\sum_{l=1}^N\left[\frac{1}{2} \boldsymbol{\mathcal{F}}(\boldsymbol{\varepsilon}(\boldsymbol{b}(\boldsymbol{\xi}_l)*\boldsymbol{\psi}(\boldsymbol{\xi}_l; \boldsymbol{\theta})))_{ij}  \boldsymbol{\varepsilon}(\boldsymbol{b}(\boldsymbol{\xi}_l)*\boldsymbol{\psi}(\boldsymbol{\xi}_l; \boldsymbol{\theta}))_{ij} - \boldsymbol{f}_0(\boldsymbol{\xi}_l) \cdot \boldsymbol{b}(\boldsymbol{\xi}_l)*\boldsymbol{\psi}(\boldsymbol{\xi}_l; \boldsymbol{\theta}) \right] \nonumber
  \\
  &- \frac{|\Gamma_T|}{N_T}\sum_{l=1}^{N_T} \big[ \boldsymbol{f}_2(\boldsymbol{\eta}_l) \cdot \boldsymbol{b}(\boldsymbol{\eta}_l)*\boldsymbol{\psi}(\boldsymbol{\eta}_l; \boldsymbol{\theta}) \big]
  + \frac{|\Gamma_C|}{N_T}\sum_{l=1}^{N_C}\big[ j_{\tau}((\boldsymbol{b}*\boldsymbol{\psi})_{\tau}(\boldsymbol{\zeta}_l; \boldsymbol{\theta})) \big],
\end{align}
where $\{\boldsymbol{\xi}_i\}_{i=1}^N$, $\{\boldsymbol{\eta}_i\}_{i=1}^{N_T}$, and $\{\boldsymbol{\zeta}_i\}_{i=1}^{N_C}$ are independent random vectors following the uniform distribution over $\Omega$, $\Gamma_T$, and $\Gamma_C$, respectively.

\vspace{0.25cm}
Similarly, for the frictionless normal compliance contact problem~\eqref{FNC3},
\begin{align}
& \quad E(\boldsymbol{b}(\boldsymbol{x})*\boldsymbol{\psi}(\boldsymbol{x}; \boldsymbol{\theta}))  \nonumber
\\
   &= \frac{1}{2}(\boldsymbol{\mathcal{F}}(\boldsymbol{\varepsilon}(\boldsymbol{b}(\boldsymbol{x})*\boldsymbol{\psi}(\boldsymbol{x}; \boldsymbol{\theta})), \boldsymbol{\varepsilon}(\boldsymbol{b}(\boldsymbol{x})*\boldsymbol{\psi}(\boldsymbol{x}; \boldsymbol{\theta}))_Q
   + \int_{\Gamma_C} j_{\nu}((b*\psi)_{\nu}(\boldsymbol{x}; \boldsymbol{\theta})) \dd s
   - \langle \boldsymbol{f}, \boldsymbol{b}(\boldsymbol{x})*\boldsymbol{\psi}(\boldsymbol{x}; \boldsymbol{\theta})) \rangle  \nonumber
   \\
   &=|\Omega|\mathbb{E}_{\boldsymbol{\xi}_1}\left[\frac{1}{2} \boldsymbol{\mathcal{F}}(\boldsymbol{\varepsilon}(\boldsymbol{b}(\boldsymbol{\xi}_1)*\boldsymbol{\psi}(\boldsymbol{\xi}_1; \boldsymbol{\theta})))_{ij}  \boldsymbol{\varepsilon}(\boldsymbol{b}(\boldsymbol{\xi}_1)*\boldsymbol{\psi}(\boldsymbol{\xi}_1; \boldsymbol{\theta}))_{ij} - \boldsymbol{f}_0(\boldsymbol{\xi}_1) \cdot \boldsymbol{b}({\boldsymbol{\xi}_1})*\boldsymbol{\psi}(\boldsymbol{\xi}_1; \boldsymbol{\theta}) \right]  \nonumber
   \\ &\quad - |\Gamma_T|\mathbb{E}_{\boldsymbol{\xi}_2} \big[ \boldsymbol{f}_2(\boldsymbol{\xi}_2) \cdot \boldsymbol{b}(\boldsymbol{\xi}_2)*\boldsymbol{\psi}(\boldsymbol{\xi}_2; \boldsymbol{\theta}) \big]
   +|\Gamma_C|\mathbb{E}_{\boldsymbol{\xi}_3}\big[ j_{\nu}((b*\psi)_{\nu}(\boldsymbol{\xi}_3; \boldsymbol{\theta})) \big],
\end{align}
with $\boldsymbol{\xi}_1$, $\boldsymbol{\xi}_2$, and $\boldsymbol{\xi}_3$ being random vectors following the uniform distribution over $\Omega$, $\Gamma_T$, and $\Gamma_C$, respectively.
In the implementation of stochastic optimization algorithm, the loss function evaluated in each iteration is defined as
\begin{align}\label{FNC5}
&\quad \hat{E}(\boldsymbol{b}(\boldsymbol{x})*\boldsymbol{\psi}(\boldsymbol{x}; \boldsymbol{\theta})) \nonumber
  \\
  &=\frac{|\Omega|}{N}\sum_{l=1}^N\left[\frac{1}{2} \boldsymbol{\mathcal{F}}(\boldsymbol{\varepsilon}(\boldsymbol{b}(\boldsymbol{\xi}_l)*\boldsymbol{\psi}(\boldsymbol{\xi}_l; \boldsymbol{\theta})))_{ij}  \boldsymbol{\varepsilon}(\boldsymbol{b}(\boldsymbol{\xi}_l)*\boldsymbol{\psi}(\boldsymbol{\xi}_l; \boldsymbol{\theta}))_{ij} - \boldsymbol{f}_0(\boldsymbol{\xi}_l) \cdot \boldsymbol{b}(\boldsymbol{\xi}_l)*\boldsymbol{\psi}(\boldsymbol{\xi}_l; \boldsymbol{\theta}) \right] \nonumber
  \\
  &- \frac{|\Gamma_T|}{N_T}\sum_{l=1}^{N_T} \big[ \boldsymbol{f}_2(\boldsymbol{\eta}_l) \cdot \boldsymbol{b}(\boldsymbol{\eta}_l)*\boldsymbol{\psi}(\boldsymbol{\eta}_l; \boldsymbol{\theta}) \big]
  + \frac{|\Gamma_C|}{N_C}\sum_{l=1}^{N_C}\big[ j_{\nu}(({b}*{\psi})_{\nu}(\boldsymbol{\zeta}_l; \boldsymbol{\theta})) \big],
\end{align}
where $\{\boldsymbol{\xi}_i\}_{i=1}^N$, $\{\boldsymbol{\eta}_i\}_{i=1}^{N_T}$, and $\{\boldsymbol{\zeta}_i\}_{i=1}^{N_C}$ are independent random vectors following the uniform distribution over $\Omega$, $\Gamma_T$, and $\Gamma_C$, respectively.

\subsection{Numerical results}
In our numerical experiments, we apply the ResNet $\boldsymbol{\phi}(\boldsymbol{x};\boldsymbol{\theta})$ or the special block ResNet $\boldsymbol{B}(\boldsymbol{x};\boldsymbol{\theta})$ introduced in Section~\ref{sec: DL} as the solution ansatz to HVIs. In the ResNet, we set the depth of $\boldsymbol{\phi}$ as $L=8$ and the width as $N=50$. In the block ResNet $\boldsymbol{B}$, we take a ResNet with depth $L=4$ and width $N=50$ as the input block and take a ResNet with depth $L=4$ and width $N=10$ as each part in the parallel block. Therefore, the total number of parameters in $\boldsymbol{\phi}(\boldsymbol{x};\boldsymbol{\theta})$ is about $20,400$ and the total number of $\boldsymbol{B}(\boldsymbol{x};\boldsymbol{\theta})$ is about $12,400$.
Note that if we choose $\sigma(x)=\text{ReLU}(x)$, the derivative of the corresponding DNN is a constant almost everywhere, which leads to DNNs not being able to capture the feature of the PDE solution. In the literature, the activation function $\text{ReLU}^{\alpha}$ or $\tanh$ are used for different problem, where $\text{ReLU}(x) = \max\{x,0 \}$ and $\alpha$ is a positive integer, to overcome this difficulty.
All neural networks are trained by Adam optimizer~\cite{KingmaBa2014} with a default learning rate $\eta = 0.001$ and exponential decay rates $\beta_1 = 0.9$ and $\beta_2 = 0.999$. The batch size in the domain is 1024 and the number of training dates on each boundary is 256 for all problems.
All numerical experiments are implemented in Python 3.7 using Pytorch 1.3 in an Nvidia GEFORCE RTX 2080 Ti GPU card.

Before reporting numerical results, let us summarize notations used in this section. Denote by $\boldsymbol{u}^{DL}=(u^{DL}_1,u^{DL}_2)^T$ as the approximate solution estimated by the deep learning method. Suppose $\mathcal{T}_h$ is a uniform triangulation of $\Omega$ into $K$ and $h=$diam$(K)$.  Since the true solution is unavailable, for the bilateral contact problem, we use $\boldsymbol{u}^{ref}$ as the reference solution evaluated by VEM with $h = 2^{-7}$ (see~\cite{FengHanHuang2019}). Similarly, for the frictionless contact with normal compliance problem, we take $\boldsymbol{u}^{ref}$ as the reference solution evaluated by VEM with $h = 2^{-9}$ (see~\cite{FengHanHuang2019}). Denote the relative difference between the deep learning solution and the reference solution as
\[
   \mathcal{E}_r = \frac{\| \boldsymbol{u}^{DL} - \boldsymbol{u}^{ref}\|_E}{\| \boldsymbol{u}^{ref}\|_E},
\]
where the energy norm $\| \cdot \|_E$ is defined by
\[
  \| \boldsymbol{v} \|_E = \frac{1}{\sqrt{2}}\left(\boldsymbol{\mathcal{F}}(\boldsymbol{\varepsilon}(\boldsymbol{v})), \boldsymbol{\varepsilon}(\boldsymbol{v}) \right)_Q^{1/2}.
\]

\vspace{0.25cm}
\begin{example}
Consider the bilateral contact problem~\eqref{FBP2}. The domain $\Omega = (0,4) \times (0,4)$ is the cross section of a three-dimensional linearly elastic body and the plane stress condition is imposed. The body is clamped on $\Gamma_D =\{ 4 \} \times (0,4)$ and vertical tractions act on $\Gamma_T =( \{ 0 \} \times (0,4) ) \cup ( (0,4) \times \{ 4 \})$. The frictional contact happens on the boundary $\Gamma_C = (0,4) \times \{ 0 \}$. The linear elasticity tensor $\boldsymbol{\mathcal{F}}$ is
\[
  \mathcal{F}_{ij} = \frac{E\kappa}{(1-\kappa^2)}(\varepsilon_{11} + \varepsilon_{22})\delta_{ij} + \frac{E}{1+\kappa}\varepsilon_{ij}, \quad 1\le i,j, \le 2,
\]
where $E$ is the Young modulus, $\kappa$ is the Poisson ratio of the material, and $\delta_{ij}$ is the Kronecker symbol. In numerical simulations, we use the following data
\begin{align*}
   &E = 2000 \ \text{daN}/\text{mm}^2, \quad \kappa = 0.4,
   \\
   &\boldsymbol{f}_0(x,y) = (0,0)^T \ \text{daN}/\text{mm}^2 \ \text{in} \ \Omega,
   \\
   &\boldsymbol{f}_2(x,y) =
   \begin{cases}
     (200(5-y),-200)^T \ \text{daN}/\text{mm}^2 \quad \text{on} \ \{ 0 \} \times (0,4)
     \\
     (0,0)^T \ \text{daN}/\text{mm}^2 \quad \text{on} \ (0,4) \times \{ 4 \},
   \end{cases}
   \\
   & j_{\tau}(\boldsymbol{z}) = \int_0^{\| \boldsymbol{z}\|} 450e^{-2000 t} + 450 \dd t, \quad \boldsymbol{z} = (x,y).
\end{align*}
\end{example}

\begin{figure}[t]
  \centering
  \includegraphics[width = 15 cm]{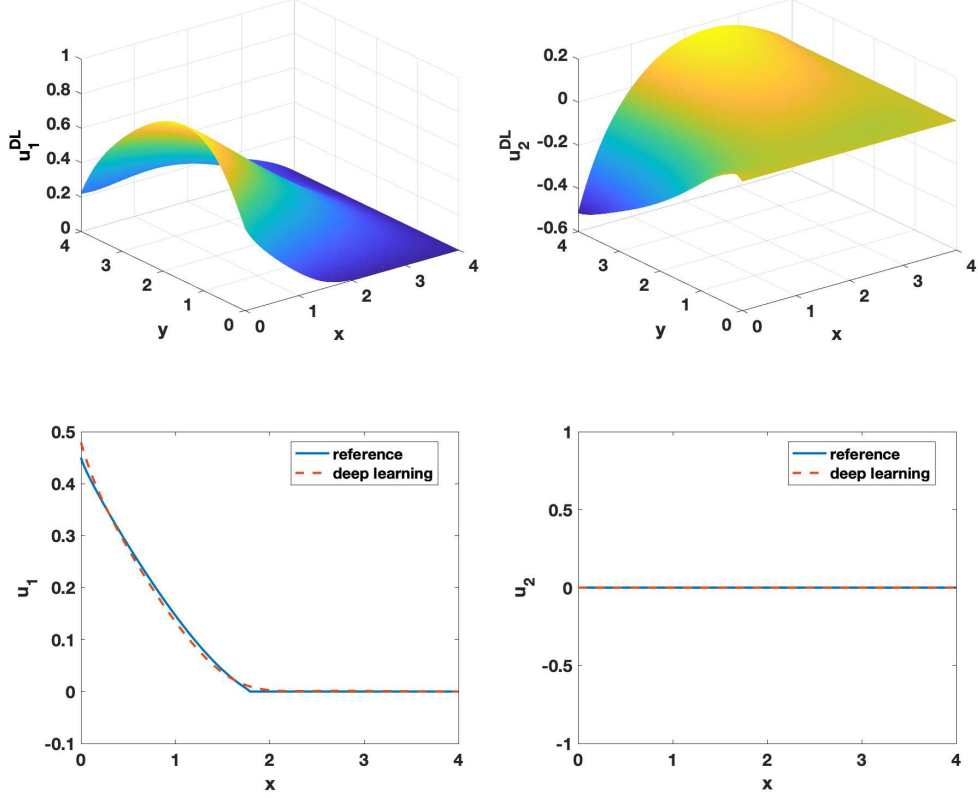}
  \caption{The numerical solution in the domain (upper) and on the contact boundary (bottom) of the bilateral contact problem.}
  \label{fig:sol1}
\end{figure}

Since that $\text{ReLU}^{\alpha}$ activation function is still not appropriate in this problem because of the numerical overflow in the training process.  We choose the tanh activation function in this example. Figure~\ref{fig:sol1} shows the numerical solution in the form of ResNet after $50,000$ epochs. The first row of Figure~\ref{fig:sol1} displays each component of $\boldsymbol{u}^{DL}$ in the domain whereas the second row shows the approximation on the contact boundary.

For the purpose of quantifying the accuracy of deep learning methods and comparing the efficiency of different networks and algorithms, we evaluate the relative error between the numerical solution and a reference solution obtained by VEM. For Algorithm~\ref{alg: basic}, we train the network with $50,000$ epochs either using the ResNet $\boldsymbol{\phi}$ or the block ResNet $\boldsymbol{B}$. For Algorithm~\ref{alg: blockwise}, we set $Epoch_{int}=9000$, $Epoch_{b}=1000$, and $Epoch_{re}=9$. In order to make fair comparisons, we stop Algorithm~\ref{alg: blockwise} after $50,000$ total epochs to update network parameters. For Algorithm~\ref{alg: adaptive}, we set the discretization step size $H=1/50$ and use $Epoch_{int}=9000$, $Epoch_{re}=41$, and $Epoch_{b}=1000$.

From Table~\ref{tab:DL error1}, we see that the relative errors are reduced to $5\%$ nearly by all training algorithms. Besides, the number of parameters of ResNet $\boldsymbol{\phi}$ is almost twice of the one of block ResNet $\boldsymbol{B}$, but the solution in the form of block ResNet $\boldsymbol{B}$ is more accurate than the one of ResNet $\boldsymbol{\phi}$. Moreover, the adaptive mesh-free multigrid algorithm can provide a  more accurate solution than others after the same number of epochs: $41.4\%$ higher accuracy than the basic algorithm and $35.5\%$ higher accuracy than the blockwise algorithm.

\begin{table}[H]
\small
\centering
\caption{The relative error of different algorithms for the bilateral contact problem. Algorithm~\ref{alg: basic}: the basic algorithm in the literature. Algorithm~\ref{alg: blockwise}: blockwise training. Algorithm~\ref{alg: adaptive}: adaptive mesh-free multigrid training. }
  \begin{tabular}{ccccc}
    \hline
    Algorithm &Algorithm~\ref{alg: basic} (ResNet $\boldsymbol{\phi}$) &Algorithm~\ref{alg: basic} (Block ResNet $\boldsymbol{B}$) &Algorithm~\ref{alg: blockwise} &Algorithm~\ref{alg: adaptive} \\
    \hline
    $\mathcal{E}_r$ &0.0481  &0.0412 &0.0437 &0.0282\\
    \hline
  \end{tabular}
  \label{tab:DL error1}
\end{table}

\vspace{0.25cm}
\begin{example}
Consider the frictionless normal compliance contact problem~\eqref{FNC2}. The domain $\Omega = (0,1) \times (0,1)$ is the cross section of a three-dimensional linearly elastic body and the plane stress condition is imposed. The body is clamped on $\Gamma_D =(\{ 0 \} \times (0,1) ) \cup (\{ 1 \} \times (0,1) )$ and the vertical traction acts on $\Gamma_T =\ (0,1) \times \{ 1 \}$. The frictional contact happens on the boundary $\Gamma_C = (0,1) \times \{ 0 \}$. The linear elasticity tensor $\boldsymbol{\mathcal{F}}$ is
\[
  \mathcal{F}_{ij} = \frac{E\kappa}{(1+\kappa)(1-2\kappa)}(\varepsilon_{11} + \varepsilon_{22})\delta_{ij} + \frac{E}{1+\kappa}\varepsilon_{ij}, \quad 1\le i,j \le 2,
\]
where $E$ is the Young modulus, $\kappa$ is the Poisson ratio of the material, and $\delta_{ij}$ is the Kronecker symbol. In numerical simulations, we use the following data
\begin{align*}
   &E = 70 \ \text{GPa}, \quad \kappa = 0.3,
   \\
   &\boldsymbol{f}_0(x,y) = (0,0)^T \ \text{GPa} \quad \text{in} \ \Omega   \\
   &\boldsymbol{f}_2(x,y) = (0,-52)^T \ \text{GPa} \quad \text{on} \ \Gamma_T \\
   & j_{\nu}(u_{\nu}) =
   \begin{cases}
     0, \quad u_{\nu} \le 0,\\
     50 u_{\nu}^2 + 0.1u_{\nu}, \quad u_{\nu} \in (0,0.1],\\
     20.1 u_{\nu} - 50 u_{\nu}^2 - 1, \quad u_{\nu} \in (0.1,0.15),\\
     200 u_{\nu}^2 - 54.9u_{\nu} + 4.625, \quad u_{\nu}\ge 0.15.
   \end{cases}
\end{align*}
\end{example}

\begin{figure}[H]
  \centering
  \includegraphics[width = 15 cm]{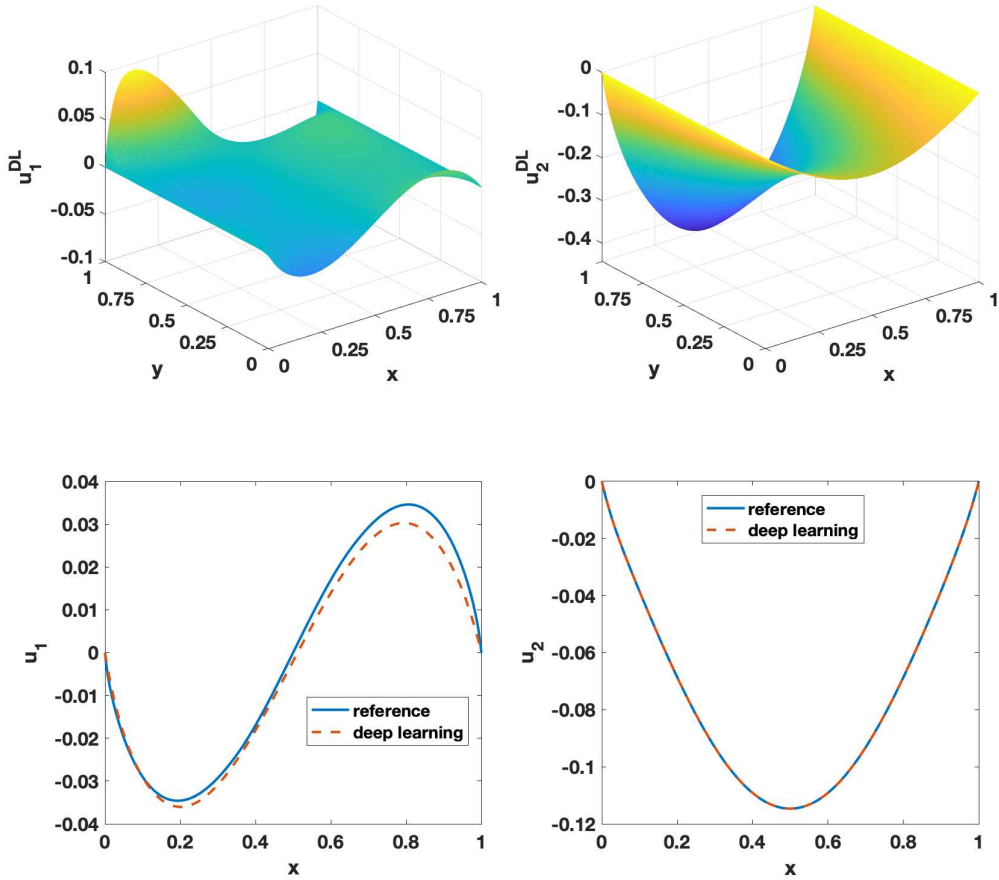}
  \caption{The numerical solution in the domain (upper) and on the contact boundary (bottom) of the frictionless normal compliance contact problem.}
    \label{fig:sol2}
\end{figure}

 The activation function $\text{ReLU}^2=\max\{x^2,0\}$ is chosen in this example based on numerical experience in this example. Figure~\ref{fig:sol2} shows the numerical solution in the form of ResNet with activation function $\sigma(x) = \text{ReLU}^2 $ after $50,000$ epochs. The first row of Figure~\ref{fig:sol2} displays each component of $\boldsymbol{u}^{DL}$ in the domain whereas the second row shows the numerical solution on the contact boundary.

Similarly, we evaluate the relative error between the numerical solution and the reference solution to compare the numerical performance of different networks and different training algorithms. For Algorithm~\ref{alg: basic}, we train the network with $50,000$ epochs either using the ResNet $\boldsymbol{\phi}$ or the block ResNet $\boldsymbol{B}$. For Algorithm~\ref{alg: blockwise}, we set $Epoch_{int}=9000$, $Epoch_{b}=1000$, and $Epoch_{re}=9$. In order to make a fair comparison, we stop Algorithm~\ref{alg: blockwise} after $50,000$ total epochs. For Algorithm~\ref{alg: adaptive}, we set the discretization step size $H=1/200$ and use $Epoch_{int}=9000$, $Epoch_{re}=41$, and $Epoch_{b}=1000$.

From Table~\ref{tab:DL error2}, we see that relative errors are reduced to $7\%$ by all deep learning methods. However, only the adaptive mesh-free multigrid algorithm can provide an approximate solution with a relative error less than $5\%$ after $50,000$ epochs. In particular, the adaptive mesh-free multigrid algorithm improves the accuracy by $38.4\%$ compared to the basic algorithm and by $22.0\%$ compared to the blockwise algorithm. At last, recall that the number of parameters of the ResNet $\boldsymbol{\phi}$ is twice the one of the block ResNet $\boldsymbol{B}$, but the approximate solution in the form of the block ResNet $\boldsymbol{B}$ is more accurate than the one of the ResNet $\boldsymbol{\phi}$.

\begin{table}[H]
\small
\centering
\caption{The relative error of different algorithms for the frictionless normal compliance contact problem. Algorithm~\ref{alg: basic}: the basic algorithm in the literature. Algorithm~\ref{alg: blockwise}: blockwise training. Algorithm~\ref{alg: adaptive}: adaptive mesh-free multigrid training. }
  \begin{tabular}{ccccc}
    \hline
    Algorithm &Algorithm~\ref{alg: basic} (ResNet $\boldsymbol{\phi}$) &Algorithm~\ref{alg: basic} (Block ResNet $\boldsymbol{B}$) &Algorithm~\ref{alg: blockwise} &Algorithm~\ref{alg: adaptive} \\
    \hline
    $\mathcal{E}_r$ &0.0706  &0.0556 &0.0558 &0.0435 \\
    \hline
  \end{tabular}
  \label{tab:DL error2}
\end{table}

\section{Conclusion}\label{sec:conclusion}
This paper focuses on developing a deep learning method for solving HVIs based on their variational problems. First of all, the solution space is parameterized via DNNs, and then the HVI is reformulated as an expectation minimization problem, which are therefore worked out by stochastic gradient descent method or its variants (e.g. Adam) combined with three different training strategies for updating network parameters. As applications to contact mechanics, a frictional bilateral contact problem and a frictionless normal compliance contact problem are carried out in details. Numerical results show that the deep learning method is efficient in solving HVIs and the adaptive mesh-free multigrid algorithm can provide the best accuracy solution among the three learning methods discussed. It deserves to mention that the proposed method is easy to realize in programming.

\section*{Acknowledgments}
J. H. was partially supported by the National Key Research and Development Project (2020YFA0709800) and NSFC (Grant No. 12071289). C. W. was partially supported by National Science Foundation Award DMS-1849483. We would like to thank Dr. Fang Feng for offering the comparison data of numerical solutions for contact problems by the virtual element method used in numerical experiments.

\bibliographystyle{plain}\footnotesize

\end{document}